\documentclass[letterpaper, 10 pt, conference]{ieeeconf}  

\IEEEoverridecommandlockouts                              

\usepackage{amsmath} 
\usepackage{amssymb}  
\usepackage{amsfonts}
\usepackage{amsthm}
\newtheoremstyle{indented}
    {3pt}
    {3pt}
    {\addtolength{\leftskip}{0em}}
    {}
    {\bfseries}
    {.}
    {1em}
    {}
\theoremstyle{indented}

\newtheorem{theorem}{Theorem}[section]
\newtheorem{proposition}{Proposition}[section]

\newtheorem{exmp}{Example}[section]
\newtheorem{remark}{Remark}[section]
\newtheorem{lemma}{Lemma}[section]

\usepackage{hyperref}
\usepackage{xcolor}
\hypersetup{
    colorlinks,
    linkcolor={blue!100!black},
    citecolor={blue!50!black},
    urlcolor={blue!80!black}
}

\newcommand{\bc}{\begin{center}}
\newcommand{\ec}{\end{center}}
\newcommand{\benum}{\begin{enumerate}}
\newcommand{\eenum}{\end{enumerate}}
\newcommand{\nn}{\nonumber}
\newcommand{\matl}{\left[ \begin{array}}
\newcommand{\matr}{\end{array} \right]}
\newcommand{\matls}{\left[ \begin{smallmatrix}}
\newcommand{\matrs}{\end{smallmatrix} \right]}
\newcommand{\isdef}{\stackrel{\triangle}{=}}

\newcommand{\rmT}{{\rm T}}

\newcommand{\rmc}{{\rm c}}

\newcommand{\rmn}{{\rm n}}

\newcommand{\rmz}{{\rm z}}

\newcommand{\BBC}{{\mathbb C}}

\newcommand{\BBR}{{\mathbb R}}

\newcommand{\real}{{\mathbb{R}}}

\newcommand{\SA}{{\mathcal A}}
\newcommand{\SB}{{\mathcal B}}
\newcommand{\SC}{{\mathcal C}}

\newcommand{\SN}{{\mathcal N}}

\newcommand{\SR}{{\mathcal R}}

\newcommand{\SZ}{{\mathcal Z}}

\renewcommand{\matl}{\begin{bmatrix}}
\renewcommand{\matr}{\end{bmatrix} }

\newcommand{\adj}{{\rm adj  } }
\newcommand{\rank}{{\rm rank  } }
\newcommand{\zeros}{{\rm zeros} }
\newcommand{\poles}{{\rm poles} }

\newcommand{\realization}[4]{
\left[\begin{array}{c|c}
    #1 & #2\\
    \hline
    #3 & #4 \end{array}\right]
}

\newcommand{\EndProofInEq}{\tag*{\mbox{\qed}}}

\usepackage[style=ieee ,sorting=none,maxbibnames=99,giveninits, doi=false]{biblatex}
\title{\LARGE \bf
Zeros in the State-Space Realization
}

\title{\LARGE \bf
Zeros are the Poles of the Zero Dynamics in Linear Systems
}

\title{\LARGE \bf
Computing Invariant Zeros of a Linear System
\\
Using State-Space Realization
}

\author{Jhon Manuel Portella Delgado and Ankit Goel
\thanks{Jhon Manuel Portella Delgado is a Ph.D. student in the Department of Mechanical Engineering, University of Maryland, Baltimore County, 1000 Hilltop Circle, Baltimore, MD 21250. 
{\href{mailto:jportel1@umbc.edu}{\tt\small jportella@umbc.edu}}}%
\thanks{Ankit Goel is an Assistant Professor in the Department of Mechanical Engineering, University of Maryland, Baltimore County,1000 Hilltop Circle, Baltimore, MD 21250. 
{\href{mailto:ankgoel@umbc.edu}{\tt\small ankgoel@umbc.edu} }}%
}

\begin{document}

\maketitle

\begin{abstract}
It is well known that zeros and poles of a single-input, single-output system in the transfer function form are the roots of the transfer function's numerator and the denominator polynomial, respectively.  
However, in the state-space form, where the poles are a subset of the eigenvalue of the dynamics matrix and thus can be computed by solving an eigenvalue problem, the computation of zeros is a non-trivial problem. 
%
%
%
This paper presents a realization of a linear system that allows the computation of invariant zeros by solving a simple eigenvalue problem. 
The result is valid for square multi-input, multi-output (MIMO) systems, is unaffected by lack of observability or controllability, and is easily extended to wide MIMO systems. 
%
Finally, the paper illuminates the connection between the zero-subspace form and the normal form to conclude that \textit{zeros are the poles of the system's zero dynamics}.

%
\end{abstract}

\section{INTRODUCTION}

Zeros are fundamental in the study of systems and control theory. 
While the poles affect the system stability, transients, and convergence rate, zeros affect the undershoot, overshoot, and zero crossings
\cite{macfarlane1976poles,desoer1974zeros,tokarzewski2006finite}.
Furthermore, nonminimum phase zeros, which are zeros in the open-right-half-plane, limit performance and bandwidth due to limited gain margin, exacerbate the tradeoff between the robustness and achievable performance of a feedback control system, and prevent input-output decoupling \cite{hoagg2007nonminimum,havre2001achievable,wonham1970decoupling}. 
Precise knowledge of zeros is thus crucial in the design of reliable control and estimation systems. 

In the transfer function representation of single-input, single-output (SISO) linear systems, 
zeros and poles are the roots of the numerator and the denominator polynomials, respectively.
In the state-space form of a SISO as well as a MIMO linear system, the system's poles are a subset of the eigenvalues of the dynamics matrix.
%
However, the zeros are not readily apparent in the state-space form. 
Invariant zeros of a MIMO system with a state-space realization $(A,B,C,D)$, which are the complex numbers for which the rank of the Rosenbrock system matrix  
\begin{align}
    \SZ (\lambda) 
        =
            \matl 
                \lambda I - A & B \\C & -D
            \matr,
\end{align}
drops, are a subset of the generalized eigenvalues of the Rosenbrock system matrix \cite{laub1978calculation,tokarzewski2011invariant,tokarzewski2009zeros}. 
%
%
%

The generalized eigenvalues of the Rosenbrock system matrix can be computed by decomposing it into the generalized Schur form as shown in \cite{emami1982computation, misra1994computation,misra1989computation,misra1989minimal}. 
However, this approach yields extraneous zeros, which are removed heuristically. 
%
%
Alternatively, since zeros are invariant under output feedback and the closed-loop poles approach the zeros of the system under high-gain output feedback, 
the zeros of $(A,B,C,D)$ are a subset of the eigenvalues of $\lim_{K \to \infty } A+ BK(I- D K)^{-1} C$ \cite{davison1978algorithm,garbow1977matrix}.
%
%
Similar to the approach based on the Schur decomposition, this approach also yields extraneous zeros. 
Furthermore, both approaches are computationally expensive.

In this paper, we present a technique to compute the invariant zeros of a MIMO system that reduces the problem to an eigenvalue problem instead of a generalized eigenvalue problem. 
Since the paper's main result does not depend on the minimality of the realization, the computed zeros are indeed the invariant zeros of the state-space realization. 
This technique is motivated by and is closely related to the normal form of a dynamic system \cite{isidori1985nonlinear,Khalil:1173048}. 
Specifically, we construct a diffeomorphism that isolates the zeros in a partition of the transformed dynamics matrix. 
Next, we show that the zeros of the system are precisely the eigenvalues of this partition.
Finally, we show that the partitioned system is in fact the zero dynamics of the system. 
This observation is stated in \cite[p.~514]{Khalil:1173048}, however, the numerator polynomial of the transfer function is used to construct the state-space realization of the zero dynamics. 
In contrast, in this paper, we construct the zero dynamics as well as compute the zeros using the system's state-space realization, without requiring its transfer function. 

A similar geometric approach, based on differential geometry, is described in \cite{basile1992controlled,morris2010invariant} to compute the invariant zeros of the system by solving an eigenvalue problem.  
In contrast, in this paper, we construct the eigenvalue problem, whose solution provides the invariant zeros, by constructing a simple state transformation matrix.
Furthermore, we present a simplified proof that does not require differential geometry concepts and is instead based on simple algebra.

Although all results presented in the paper are applicable to square MIMO systems as well, due to page limits, we restrict the proofs to the case of SISO systems. 
The technique to compute the zeros can be easily extended to wide MIMO systems as shown in an example in Section \ref{sec:examples}.
However, the technique does not work in the case of tall MIMO systems since the construction of the zero-subspace form yields a singular transformation matrix. 
Ad-hoc techniques such as modifying the rows of the transformation matrix to make it nonsingular yield mixed and inconclusive results. 
The paper thus does not consider the case of tall MIMO systems.

The paper is organized as follows. 
Section \ref{sec:notation} presents the notation used in this paper, 
Section \ref{sec:ZSF} introduces the zero-subspace form of a strictly proper linear system, 
Section \ref{sec:main_result} presents and proves the main result of the paper, 
Section \ref{sec:normal_form} shows the connection of the zero-subspace form with the normal form, 
and 
Section \ref{sec:examples} presents numerical examples to confirm the main result of this paper.
Finally, the paper concludes with a discussion in Section \ref{sec:conclusion}.

\section{Notation}
\label{sec:notation}
Let $A \in \BBR^{n\times m}.$
Then, $A_{[i,j]}$ is the matrix obtained by removing the $i$th row and $j$th column of $A.$
Note that, if $i=0,$ then only the $j$th column is removed.
Similarly, if $j=0,$ then only the $i$th row is removed.
That is, $A_{[0,0]} = A.$
The set of integers between $n$ and $m$, where $n\leq m,$ that is, $\{n, n+1, \ldots, m\},$ is denoted by $\{ \iota \}_{n}^{m} .$ 
$0_{n \times m}$ denotes the $n \times m$ zero matrix and 
$I_n$ denotes the $n\times n$ identity matrix.




\section{zero-subspace form}
\label{sec:ZSF}
This section introduces the \textit{zero-subspace form} of a linear system.
%
The zero-subspace form is a realization in which the zeros of the system are the eigenvalues of a partition of the dynamics matrix.
%
%

%

%

Consider a linear system
\begin{align}
    \lambda x &= A x + Bu, \label{eq:state_x} \\
    y &= Cx, \label{eq:output_y}
\end{align}
where $x\in \BBR^{l_x}$ is the state, 
$u \in \BBR^{l_u}$ is the input, 
$y \in \BBR^{l_y}$ is the output, and 
the operator $\lambda$ is the time-derivative operator or the forward-shift operator. 
The relative degree of $y_i$ is denoted by $\rho_i,$ and the relative degree of $y,$ defined as $\sum_i^{l_y} \rho_i,$ is denoted by $\rho.$
%
%
%

Note that \eqref{eq:state_x}-\eqref{eq:output_y} is a strictly proper system since $D = 0.$ 
Although the main result presented in this paper requires $\rho>0,$ that is, direct feedthrough term $D=0,$ it is not a restrictive condition since, without affecting the zeros of the system, the dynamic extension of the system with an additional pole renders the direct feedthrough term $D$ of the augmented system zero. 
Example \ref{exmp:D_nonzero} considers the case of an exactly proper system.

Due to substantially lengthier and more complex proofs in the case of  multi-input, multi-output (MIMO) linear systems, we restrict our attention to single-input, single-output (SISO) linear systems in the following. 
%
%
However, the main result of this paper and the procedure to transform a realization of a MIMO system into its zero-subspace form is the same as that of a SISO system.
%

First, we construct a diffeomorphism to transform a realization into the zero-subspace form, as shown below.  
Let $\overline B \in \BBR^{l_x -1 \times l_x}$ be a full rank matrix such that $\overline B B = 0.$
Note that $B \in \SN(\overline{B}),$ and thus, in MATLAB, $\overline B$ can be computed using \texttt{null(B')}.
Define 
\begin{align}
    \overline{C}
        \isdef 
            \matl
                C\\
                CA\\
                \vdots\\
                CA^{\rho-1}
            \matr \in \BBR^{\rho \times l_x}, 
    \label{eq:Cbar_def}
\end{align}
and 
\begin{align}
    \overline T 
        \isdef
            \matl
                \overline{B} \\
                \overline{C}
            \matr
            \in \BBR^{(l_x-1+\rho) \times l_x}.
    \label{eq:Tbar_def}
\end{align}
\begin{proposition}
    \label{prop:Tbar_full_rank}
    Let $\overline T$ be given by \eqref{eq:Tbar_def}.
    Then, $\rank \ \overline T = l_x.$
\end{proposition}
\textbf{Proof.}
    Note that $\overline B$ has $l_x-1$ linearly independent rows. 
    Since the relative degree of $y$ is $\rho,$ it follows that, for 
    $ i \in \{ \iota \}_{\iota =1}^{\rho-1}, $
    $C A^{i-1} B = 0,$ which implies that each element of $\{ C, CA, \ldots, CA^{\rho-2} \}$ is in the row range space of $\overline B.$
    Furthermore, since $ C A^{\rho-1} B \neq 0,$ it follows that $C A^{\rho-1}$ is not linearly dependent on the rows of $\overline B,$ thus implying that $\rank \ \overline T = l_x.$
    \qed

Next, define $l_z \isdef l_x - \rho$ and  
\begin{align}
    T 
        \isdef 
            \matl 
                B_\rmz \\
                \overline C
            \matr 
            \in \BBR^{l_x \times l_x},
    \label{eq:T_def}
\end{align}
where 
$B_\rmz \in \BBR^{l_z \times l_x}$ is chosen such that rows of $B_\rmz$ and rows of $\overline C$ are linearly independent. 
Consequently, $ T$ is full-rank and thus invertible. 
Note that $\SR(B_\rmz^\rmT) \subseteq \SR(\overline B^\rmT).$ 
%
%
    Using the diffeomorphism $T,$ transform the realization $(A,B,C)$ into the zero-subspace form $(\SA,\SB,\SC),$ which implies that  
    \begin{align}
        \SA \isdef T A T^{-1}, \quad 
        \SB 
            \isdef
                TB, \quad 
        \SC 
            \isdef 
                C T^{-1}.
        \label{eq:ZSF_SS}
    \end{align}

    Next, define $\eta \in \BBR^{l_z}$ as the first $l_z$ components of $Tx$ and 
    $\xi \in \BBR^{\rho}$ as the rest of $Tx, $ that is, 
    \begin{align}
        \matl 
            \eta \\
            \xi 
        \matr 
            =
                T x.
        \label{eq:state_transformation}
    \end{align}
    Note that since $\eta \in \SR(B_\rmz),$ we call $\SR(B_\rmz)$ the \textit{zero subspace} of the system. 



Next, partition $\SA$ as
\begin{align}
    \SA 
        =
            \matl 
                \SA_{\eta} & \SA_{\eta \xi} \\
                \SA_{\xi \eta} & \SA_{\xi}
            \matr,
\end{align}
where $\SA_\eta$ is the $l_z \times l_z$ upper-left block,
$\SA_{\eta \xi}$ is the $l_z \times \rho$ upper-right block, 
$\SA_{\xi \eta}$ is the $\rho \times l_z$ lower-left block, and 
$\SA_\xi$ is the $\rho \times \rho$ lower-right block of $\SA.$

Finally, define $S \isdef T^{-1}$ and partition $S$ as $S = \matl S_\eta & S_\xi \matr,$ where $S_\eta \in \BBR^{l_x \times l_z}$ contains the first $l_z$ columns of $S$ and $S_\xi \in \BBR^{l_x \times \rho}$ contains the last $\rho$ columns of $S.$
Substituting $T$ and $S$ in \eqref{eq:ZSF_SS} yields
\begin{align}
    \SA_\eta  &=  B_\rmz A S_\eta,
    \label{eq:SA_eta_def}
    \\ 
    \SA_{\eta \xi}  &= B_\rmz A S_\xi, \\ 
    \SA_{\xi \eta}  &= \overline C A S_\eta, \label{eq:SA_xi_eta_def} \\ 
    \SA_\xi  &= \overline C A S_\xi. \label{eq:SA_xi_def} 
\end{align}

As will be shown in the next section, the zeros of \eqref{eq:state_x}-\eqref{eq:output_y} are the 
eigenvalues of $\SA_\eta$ given by \eqref{eq:SA_eta_def}.

\subsection{Sparse Structure of the zero-subspace form}
The following two facts about the zero-subspace form show that the matrices $\SA_{\xi\eta},$ $\SA_\xi,$ $\SB,$ and $\SC$ have \textit{sparse} structure. 

\begin{proposition}
    \label{prop:SA_xi_eta_form}
    Let $\SA_{\xi \eta}$ and $\SA_\xi$ be defined by \eqref{eq:SA_xi_eta_def} and \eqref{eq:SA_xi_def}.
    Then, 
    \begin{align}        
        A_{\xi \eta}
            &=
                \matl
                    0_{(\rho-1) \times l_z}\\
                    CA^{\rho} S_\eta
                \matr,
        \label{eq:SA_xi_eta_def_2}
    \\ 
        A_{\xi }
            &=
                \matl
                    \matl 
                        0_{(\rho-1) \times 1} & I_{\rho-1}
                    \matr \\
                    CA^{\rho} S_\xi
                \matr. 
    \label{eq:SA_xi_def_2}
    \end{align}
\end{proposition}
\textbf{Proof.}
Since $\matl B_\rmz \\ \overline C \matr \matl S_\eta & S_\xi \matr = I_{l_x},$ it follows that 
%
\begin{align}
    \overline C S_\eta         
        &=
            \matl 
                C \\
                \vdots\\
                C A^{\rho-1} 
            \matr 
            S_\eta  
        =
            0_{\rho \times l_z}, \nn
    \\
    \overline C S_\xi         
        &=
            \matl 
                C \\
                \vdots\\
                C A^{\rho-1} 
            \matr 
            S_\xi  
        =
            I_{\rho}, \nn
\end{align}
which implies that, for each 
$ j \in \{ \iota \}_{\iota =0}^{\rho-1}, $ $C  A^{j} S_\eta = 0$
and 
$CA^j S_\xi = e_{j+1},$ where $e_j$ is the $j$th row of $I_\rho.$ 
Therefore, 
\begin{align*}
    \SA_{\xi \eta}  
        &=
            \overline C A S_\eta
        = 
            \matl 
                CA S_\eta\\
                \vdots\\
                C A^{\rho-1} S_\eta \\
                C A^{\rho} S_\eta
            \matr             
        =
            \matl
                0_{(\rho-1) \times l_z}\\
                C A^{\rho} S_\eta
            \matr,
        \\
    \SA_{\xi}  
        &=
            \overline C A S_\xi
        = 
            \matl 
                CA S_\xi\\
                \vdots\\
                C A^{\rho-1} S_\xi \\
                C A^{\rho} S_\xi
            \matr             
        =
            \matl
                \matl 
                    0_{(\rho-1) \times 1} & I_{\rho-1}
                \matr \\
                CA^{\rho} S_\xi
            \matr.
    \EndProofInEq
\end{align*}

\begin{proposition}
    Let $\SB$ and $\SC$ be defined by \eqref{eq:ZSF_SS}.
    Then, 
    \begin{align}
    \SB 
        =
            \matl 
                0_{(l_x-1)\times 1} \\
                C A^{\rho-1}B 
            \matr, \quad 
    \SC 
        =
            \matl 
                0_{1\times l_z} &
                1
                &
                0_{1\times (\rho-1)}
            \matr.
\end{align}
\end{proposition}
\textbf{Proof.}
Note that, for each 
$ j \in \{ \iota \}_{\iota =0}^{\rho-1}, $
$C A^j B = 0,$ which implies 
\begin{align}
    \SB
        =
            TB
        =
            \matl
                B_\rmz B\\
                \overline C B
            \matr 
        =
            \matl
                0_{l_z \times l_u}\\
                \matl
                    C\\
                    CA\\
                    \vdots\\
                    CA^{\rho-1}
                \matr
                B
            \matr 
        =
            \matl
                0_{(l_x-1) \times 1}\\
                CA^{\rho-1} B
            \matr.
            \nn
\end{align}
Next, since $C$ is the $(l_z+1)$th row of $T,$ and
\begin{align*}
    C 
        =
            \SC T 
        =
            \SC  
            \matl
                B_\rmz \\
                \overline C
            \matr 
        =
            \SC
            \matl
                B_\rmz \\
                C\\
                CA\\
                \vdots\\
                CA^{\rho-1}
            \matr, 
\end{align*}
it follows that $\SC = \matl
                0_{1 \times l_z}
                &
                1
                & 
                \cdots 
                & 
                0_{1 \times \rho -1}
            \matr .$ \qed

In summary, the zero-subspace form of \eqref{eq:state_x}-\eqref{eq:output_y} is
\begin{align}
    \dot \eta 
        &=
            A_\eta \eta + A_{\eta \xi} \xi,
    \label{eq:ZSF_eta_dot}
    \\
    \dot \xi 
        &=
            A_{\xi \eta} \eta + A_\xi \xi + B_\xi u, 
    \label{eq:ZSF_xi_dot}
    \\
    y 
        &=
            C_\xi \xi = \xi_1,
    \label{eq:ZSF_y}
\end{align}
where 
\begin{align}
    B_\xi 
        \isdef 
            \matl 
                0_{(\rho-1)\times 1} \\
                C A^{\rho-1}B 
            \matr
            \in \BBR^{\rho}
            ,
    \quad 
    C_\xi 
        \isdef 
            \matl 
                1
                &
                0_{1\times (\rho-1)}
            \matr
            \in \BBR^{1\times \rho}.
\end{align}

\section{Main Result}
\label{sec:main_result}
This section presents the paper's main result, that is, the zeros of a SISO linear system are the eigenvalues of a partition of the dynamics matrix represented in its zero-subspace form.
Specifically, Theorem \ref{theorem:zeros} shows that the zeros of the system \eqref{eq:state_x}-\eqref{eq:output_y} are exactly the eigenvalues of $A_\eta.$

The following lemma appears in \cite{powell2011calculating} and is used in the proof of Theorem \ref{theorem:zeros}.
\begin{lemma}
    Let $A = \matl A_{11} & A_{12} \\ A_{21} & A_{22} \matr,$ where $A_{11}$ is square and $A_{22}$ is nonsingular. 
    Then, 
    \begin{align}
        \det A
            =
                \det 
                \left( 
                    A_{11}  - A_{12} A_{22}^{-1} A_{21}
                \right)
                \det A_{22}.
    \end{align}
\end{lemma}

The following two facts about the zero-subspace form are 
 used in the proof of Theorem \ref{theorem:zeros}.

\begin{proposition}
    \label{prop:SA_xi_eta_rho0}
    Let $\SA_{\xi \eta}$ be defined by \eqref{eq:SA_xi_eta_def}.
    Then, 
    \begin{align}
        {\SA}_{\xi \eta_{[\rho,0]}} = 0_{(\rho-1) \times l_z}.
        \label{eq:SA_xi_eta_rho0}
    \end{align}
\end{proposition}
\textbf{Proof.}
    Note that $\SA_{\xi \eta}$ is given by \eqref{eq:SA_xi_eta_def_2}.
    Removing the $\rho$th row of ${\SA}_{\xi \eta}$ yields \eqref{eq:SA_xi_eta_rho0}.
    {\hfill {\qed}}

\begin{proposition}
    \label{prop:sI-Axi_rho1}
    Let $\SA_\xi$ be defined by \eqref{eq:SA_xi_def}.
    Then, for all $s\in \BBC,$ 
    \begin{align}
        \det \left( (sI_\rho - \SA_\xi)_{[\rho,1]} \right)= \pm 1.
        \label{eq:sI-Axi_rho1}
    \end{align}
\end{proposition}
\textbf{Proof.}
Note that $\SA_{\xi}$ is given by \eqref{eq:SA_xi_def_2}.
Removing the $\rho$th row and the first column of $\SA_\xi$ yields
\begin{align}
    (sI_\rho - \SA_\xi)_{[\rho,1]}
        &=
            \matl
                -1 & 0 & \cdots & 0 & 0\\
                s & -1 & \cdots & 0 & 0\\
                \vdots & \vdots & \ddots & \vdots & \vdots\\
                0 & 0 & \cdots & -1 & 0\\
                0 & 0 & \cdots & s & -1
            \matr, \nn
\end{align}
which implies \eqref{eq:sI-Axi_rho1}.
{\hfill {\qed}}

\begin{theorem}
\label{theorem:zeros}
Consider the system \eqref{eq:state_x}-\eqref{eq:output_y}
and its zero-subspace form \eqref{eq:ZSF_eta_dot}-\eqref{eq:ZSF_y}.
Then, the zeros of the system \eqref{eq:state_x}-\eqref{eq:output_y} are the eigenvalues of $\SA_\eta$. 
\end{theorem}
\textbf{Proof.}
Note that the zeros of the system \eqref{eq:state_x}-\eqref{eq:output_y} satisfy
\begin{align}
   C \adj (sI-A) B=0. \nn 
\end{align}
Since the zeros are invariant under state transformation, it follows that the zeros also satisfies
\begin{align}
    \SC \adj (sI- \SA ) \SB=0. 
    \label{eq:zeros_adj_equation}
\end{align}
Next, note that 
\begin{align}
    \SC {\rm adj  }&(sI- \SA ) \SB\nn\\
        =&
            \matl 
                0_{1\times l_z} &
                1
                &
                0_{1\times (\rho-1)}
            \matr.\nn\\
            &\adj
            \left(
            \matl 
                sI_{l_z}-\SA_{\eta} & -\SA_{\eta \xi} \\
                -\SA_{\xi \eta} & s I_{\rho} -\SA_{\xi}
            \matr 
            \right) 
            \matl 
                0_{(l_x-1)\times 1} \\ CA^{\rho-1}B
            \matr
            \nn
            \\
        =&
            (-1)^n.\nn\\
            &\det \left(
            \matl
                sI_{l_z}-\SA_{\eta} & 
                -
                {\SA}_{{\eta \xi}_{[0,1]}}
                \\
                -
                {\SA}_{\xi \eta_{[\rho,0]}}
                &
                (sI_\rho - \SA_\xi)_{[\rho,1]}
            \matr
            \right)
            CA^{\rho-1}B
        \nn \\
        =& 
            (-1)^n
        \det
            \left(
                sI_{l_z}-\SA_{\eta} -
                {\SA}_{{\eta \xi}_{[0,1]}}
                (sI_\rho - \SA_\xi)_{[\rho,1]}^{-1}
                {\SA}_{\xi \eta_{[\rho,0]}}
            \right)
            \cdot 
        \nn \\ &\quad \quad 
        \det
            (
                (sI_\rho - \SA_\xi)_{[\rho,1]}
            )
            CA^{\rho-1}B, \nn 
        \nn \\
        =& 
            (-1)^n
        \det
            \left(
                sI_{l_z}-\SA_{\eta}
            \right)
        \det
            (
                (sI_\rho - \SA_\xi)_{[\rho,1]}
            )
            CA^{\rho-1}B, \nn 
\end{align}
where 
$n \isdef {2l_x-\rho+1}.$
Note that the last equality uses the fact that ${\SA}_{\xi \eta_{[\rho,0]}} = 0,$ which follows from Proposition \ref{prop:SA_xi_eta_rho0}.
%
%
Since the relative degree of the system is $\rho,$ $CA^{\rho-1}B \neq 0.$
Next, it follows from Proposition \ref{prop:sI-Axi_rho1} that, for all $s\in \BBC,$ $\det
            (
                (sI_\rho - \SA_\xi)_{[\rho,1]}
            ) 
            =
                \pm 1 
            \neq
                0.$
Therefore, \eqref{eq:zeros_adj_equation} is satisfied if and only if $\det (sI_{l_z}-\SA_{\eta} ) = 0,$ thus implying that the zeros of the system \eqref{eq:state_x}-\eqref{eq:output_y} are the eigenvalues of $\SA_{\eta}.$
\qed

\begin{remark}
    Note that $\SA_\eta$ has $l_z$ eigenvalues and thus \eqref{eq:state_x}-\eqref{eq:output_y} has $l_z=l_x-\rho$ zeros, which is a well-known fact.

\end{remark}

\begin{remark}
A simplistic proof of the zeros being the eigenvalues of a partition of the zero-subspace form can be motivated by the blocking property of zeros. 
%
%
%
Let $z$ denote a zero of an asymptotically stable system. 
%
It follows that the input $u(t) = e^{zt} u_0$ yields zero asymptotic output, that is, 
$
    \lim_{t \to \infty} y(t)= 0,
$
which implies that  
$
    \lim_{t \to \infty} \xi(t)= 0.
$
It follows from \eqref{eq:ZSF_xi_dot} that
\begin{align}
    \lim_{t \to \infty} \SA_{\xi \eta} \eta (t) + \SB_\xi u(t) = 0.
    \label{eq:eta_u_diverging}
\end{align}
Now, consider the case where $\real(z)>0.$
Since $u(t) = e^{zt} u_0$ diverges, it follows that $\eta(t)$ must diverge to satisfy \eqref{eq:eta_u_diverging}.
In fact, $\eta(t)$ must grow exponentially as $e^{zt}$ to satisfy \eqref{eq:eta_u_diverging}.
Since $\xi(t) \to 0,$ it follows from \eqref{eq:ZSF_eta_dot} that $\eta(t)$ grows exponentially if and only if $z \in \text{spec}(A_\eta).$


\end{remark}

\begin{theorem}
\label{theorem:zeros_MIMO}
Consider the square MIMO system \eqref{eq:state_x}, \eqref{eq:output_y}
and its zero-subspace form, given by \eqref{eq:ZSF_eta_dot}-\eqref{eq:ZSF_y}.
Then, the invariant zeros of the system \eqref{eq:state_x}, \eqref{eq:output_y} are the eigenvalues of $\SA_\eta$, that is, 
\begin{align}
        {\rm izeros } \left( \realization{A}{B}{C}{0} \right) = {\rm mspec } (\SA_\eta)\nn, 
\end{align}
where ${\rm mspec}(\SA_\eta)$ denotes the multispectrum of $\SA_\eta,$ that is, the multiset consisting of the eigenvalues of $\SA_\eta$ including their algebraic multiplicity \cite[p.~506]{bernstein2018scalar}. 

\end{theorem}

\textbf{Proof.}
The proof is omitted due to restrictions on allowable number of pages. 

\section{Normal Form of a Linear System}
\label{sec:normal_form}
This section shows the relation between the zero-subspace form \eqref{eq:ZSF_eta_dot}-\eqref{eq:ZSF_y} of a system and its normal form.
In the nonlinear systems theory, the normal form is used to decompose a nonlinear system into its zero dynamics and the input-output feedback linearizable dynamics \cite{isidori1985nonlinear,Khalil:1173048}. 
%
This section shows that a system's normal form can also be used to deduce its zeros.
%

%


%
%
To maintain consistent notation, the procedure to transform a multi-input, multi-output (MIMO), nonlinear affine system to its normal form is summarized in Appendix \hyperref[appndx:normalForm]{A}.
In the case of a SISO linear system,
$
    f(x) = Ax, \
    g(x) = B, $ and $
    h(x) = Cx,
$
and thus it follows from \eqref{eq:alpha_def}, \eqref{eq:beta_def} that
\begin{align}
    \alpha(x) 
        &=
            L_f^\rho h(x)
        =
            CA^{\rho} x,
        \\
    \beta(x) 
        &=
            L_g L_f^{\rho-1} h(x)
        =
            CA^{\rho-1} B.
\end{align}
Furthermore, 
\begin{align}
    \psi(x)
        =
            \matl 
                C \\
                CA \\
                \vdots \\
                CA^{\rho-1}
            \matr x
        . 
\end{align}
Note that $\psi(x)  = \overline C x,$ where $\overline C$ is defined by \eqref{eq:Cbar_def}.
Finally, let $\phi(x) = B_\rmn x,$ where $B_\rmn$ is chosen to satisfy $L_g \phi(x) = B_\rmn B = 0.$ 
Note that 
\begin{align}
        \matl 
            \eta \\
            \xi 
        \matr 
            =
                \matl 
                    \phi(x) \\
                    \psi(x) 
                \matr
            =
                \matl 
                    B_\rmn \\
                    \overline C
                \matr x,
        \label{eq:state_transformation_NF}
    \end{align}
and thus $x = R_\eta \eta + R_\xi \xi, $ where 
\begin{align}
    \matl R_\eta & R_\xi \matr
        =
            \matl 
                B_\rmn \\
                \overline C
            \matr^{-1}.
\end{align}
Note that the state-transformation matrix in \eqref{eq:state_transformation_NF} is invertible.
The proof is similar to the proof of Proposition \ref{prop:Tbar_full_rank}.
The normal form of \eqref{eq:state_x}-\eqref{eq:output_y} is then
\begin{align}
    \dot \eta 
        &=
            B_\rmn A x
        =
            B_\rmn A R_\eta \eta + B_\rmn A R_\xi \xi
        ,
    \label{eq:NF_eta_dot}
        \\
    \dot \xi
        &=
            B_\rmc
            CA^{\rho} 
            A R_\eta \eta 
            + 
            (A_\rmc 
            +
            B_\rmc
            CA^{\rho} 
            A R_\xi) \xi
            +
            B_\rmc
            CA^{\rho-1} B u,
    \label{eq:NF_xi_dot}
    \\
    y
        &=
            \xi_1,
\end{align}
where 
\begin{align}
    A_{\rmc}
        &\isdef 
            \matl 
                0 & 1 & 0 & \cdots & 0  \\
                0 & 0 & 1 & \cdots & 0  \\
                \vdots & \vdots & \ddots & \ddots & \vdots \\
                0 & \vdots & \ldots & 0 & 1  \\
                0 & \vdots & \ldots & 0 & 0  \\
            \matr
        \in \BBR^{\rho \times \rho},
    \quad 
    B_{\rmc}
        \isdef 
            \matl 
                0\\
                \vdots \\
                1
            \matr
        \in \BBR^{\rho}.
\end{align}

Note that since $B_\rmz$ and $B_\rmn$ satisfy exactly the same conditions, it follows that  $\SR(B_\rmz) = \SR(B_\rmn)$ and therefore \eqref{eq:ZSF_eta_dot}, \eqref{eq:ZSF_xi_dot} are same as \eqref{eq:NF_eta_dot}, \eqref{eq:NF_xi_dot}.
In fact, if $B_\rmn$ is chosen to be equal to $B_\rmz,$ then the zero-subspace form \eqref{eq:ZSF_eta_dot}, \eqref{eq:ZSF_xi_dot} is exactly the normal form of \eqref{eq:state_x}-\eqref{eq:output_y}.
%

In the normal form, \eqref{eq:NF_eta_dot} is called the \textit{zero dynamics} of the system.
Theorem \ref{theorem:zeros} thus shows that \textit{zeros are the eigenvalues of the dynamics matrix of the zero dynamics.}
\section{Numerical Examples}
\label{sec:examples}

This section presents several examples verifying the main result of this paper. 
These examples, listed in Table \ref{tab:examples}, show that the technique presented in the paper to compute the zeros can be extended to exactly proper systems and nonsquare MIMO systems and is not affected by pole-zero cancellations. 

\begin{table}[ht]
    \centering
    \renewcommand{\arraystretch}{1.5}
    \begin{tabular}{|c|l|}
        \hline
        Example & \hspace{6em} Remark  
        \\ \hline
            Example \ref{exmp:SISOsystem} &
            Strictly proper SISO system
        \\ \hline
            Example \ref{exmp:SISOsystem_w_zero_rd} & 
            Exactly proper SISO System
        \\ \hline
            Example \ref{exmp:SISOsystem_w_pzcancellation} &
            SISO system with pole-zero cancellation 
        \\ \hline
            Example \ref{exmp:square_MIMO_2_zeros} &
            Square MIMO system
        \\ \hline
            Example \ref{exmp:wide_MIMO} &
            Wide MIMO system
        \\ \hline

    \end{tabular}
    \caption{List of Examples in this paper.}
    \label{tab:examples}
\end{table}

\begin{exmp}
    \label{exmp:SISOsystem}
    \textbf{[Strictly proper SISO system.]}
    Consider the system
    \begin{align}
        G(s)=\dfrac{s^2-9s+8}{s^3+11s^2+36s+36}. \nn
    \end{align}
    Note that   $\rho=1$ and $\zeros(G) = \{ 1,8\}$.
    A controllable canonical realization of $G$ is
    \begin{align*}
        A
            &=
                \matl
                    0 & 1 & 0\\
                    0 & 0 & 1\\
                    -36 & -36 & -11\\
                \matr, \quad
        B
            =
                \matl
                    0\\
                    0\\
                    1
                \matr, 
        \\
        C
            &=
                \matl
                    8 & -9 & 1
                \matr.
    \end{align*}

    Letting $B_\rmz\footnote{Computed using \texttt{null(B')'} in MATLAB.} = \matl 0 & 1 & 0\\
    1 & 0 & 0\matr,$ it follows from \eqref{eq:T_def} that 
    \begin{align}
        T
            =
                \matl
                    B_\rmz\\
                    C
                \matr
            =
                \matl
                    0 & 1 & 0\\
                    1 & 0 & 0\\
                    8 & -9 & 1
                \matr,\nn
    \end{align}
    and thus
    \begin{align*}
        \SA
            &=
                \left[\begin{array}{cc|c}
                    9 & -8 & 1\\
                    1 & 0 & 0\\
                    \hline
                    -208 & 124 & -20
                \end{array}\right]
        \SB=\left[\begin{array}{c}
            0\\
            0\\
            \hline
            1
        \end{array}\right], 
        \\
        \SC
            &=
                \left[\begin{array}{cc|c}
                    0 & 0 & 1
                \end{array}\right].
    \end{align*}
    Note that $\SA_\eta = \matl9 & -8\\
    1 & 0\matr$ and thus ${\rm mspec}(\SA_\eta) = \{1,8\},$ which confirms Theorem \ref{theorem:zeros}.
    \hfill{\huge$\diamond$}
\end{exmp}

\begin{exmp}
\label{exmp:SISOsystem_w_zero_rd}
{\textbf{[SISO System with zero relative degree.]}}
\label{exmp:D_nonzero}
Consider the system
    \begin{align}
        G(s)=\dfrac{s^3+21s^2+116s+96}{s^3+11s^2+38s+40}.\nn
    \end{align}
    Note that   $\rho=0$ and $\zeros (G) = \{ -1,-8,-12 \}$.
    Furthermore, $l_z = 3$ and thus $\SA_\eta$ is $3 \times 3,$ which implies that $\SA_\eta$ and the dynamics matrix of $G$ are similar, that is, they have same eigenvalues.
    Therefore, in the case where $\rho=0,$ eigenvalues of $\SA_\eta$ are the poles of $G.$
    Note that Theorem \ref{theorem:zeros} is not applicable in this case since $D \neq 0.$

    However, consider the system 
    \begin{align}
        \overline G(s)
            \isdef 
                \dfrac{G(s)}{s}
            =
                \dfrac{s^3+21s^2+116s+96}{s^4+11s^3+38s^2+40s}.\nn
    \end{align}
    Note that $\overline G(s)$ has the same set of zeros as $G(s),$ but its relative degree $\rho=1$.
    Therefore, Theorem \ref{theorem:zeros} can be used to compute the zeros of $\overline G,$ and thus the zeros of $G.$
    %
    %
    %
    A controllable canonical realization of $\overline G$ is
    \begin{align*}
        A
            &=
                \matl
                    0 & 1 & 0 & 0\\
                    0 & 0 & 1 & 0\\
                    0 & 0 & 0 & 1\\
                    0 & -40 & -38 & -11\\
                \matr, \quad
        B
            =
                \matl
                    0\\
                    0\\
                    0\\
                    1
                \matr, 
        \\
        C
            &=
                \matl
                    96 & 116 & 21 & 1
                \matr.
    \end{align*}
    
    Letting $B_\rmz = \matl 0 & 1 & 0 & 0\\
    0 & 0 & 1 & 0\\
    -1 & 0 & 0 & 0\matr,$ it follows from \eqref{eq:T_def} that 
    \begin{align}
        T
            =
                \matl
                    B_\rmz\\
                    C\\
                \matr
            =
                \matl
                    0 & 1 & 0 & 0\\
                    0 & 0 & 1 & 0\\
                    -1 & 0 & 0 & 0\\
                    96 & 116 & 21 & 1
                \matr, \nn
    \end{align}
    and thus
    \begin{align}
        \SA
            &=
                \left[\begin{array}{ccc|c}
                    0 & 1 & 0 & 0\\
                    -116 & -21 & 96 & 1\\
                    -1 & 0 & 0 & 0\\
                    \hline
                    -1104 & -132 & 960 & 10
                \end{array}\right]
        \SB=\left[\begin{array}{c}
            0\\
            0\\
            0\\
            \hline
            1
        \end{array}\right], \nn
        \\
        \SC
            &=
                \left[\begin{array}{ccc|c}
                    0 & 0 & 0 & 1
                \end{array}\right].\nn
    \end{align}
    Note that $\SA_\eta = \matl
        0 & 1 & 0 \\
        -116 & -21 & 96\\
        -1 & 0 & 0
    \matr$ and thus ${\rm mspec}(\SA_\eta) = \{-12,-8,-1\},$ which confirms Theorem \ref{theorem:zeros}.
    \hfill{\huge$\diamond$}
    
\end{exmp}

\begin{exmp}
\label{exmp:SISOsystem_w_pzcancellation}
{\textbf{[SISO System with pole-zero cancellation.]}}
Consider the system
    \begin{align}
        G(s)=\dfrac{s+5}{s^3+10s^2+31s+30}.\nn
    \end{align}
    Note that   $\rho=2,$  $\zeros (G) = \{ -5 \},$ and $\poles (G) = \{-2, -3, -5\}.$
    A controllable canonical realization of $G$ is
    \begin{align*}
        A
            &=
                \matl
                    0 & 1 & 0\\
                    0 & 0 & 1\\
                    -30 & -31 & -10\\
                \matr, \quad
        B
            =
                \matl
                    0\\
                    0\\
                    1
                \matr, 
        \\
        C
            &=
                \matl
                    5 & 1 & 0
                \matr.
    \end{align*}

    Letting $B_\rmz = \matl 0 & 1 & 0 \matr,$ it follows from \eqref{eq:T_def} that 
    \begin{align}
        T
            =
                \matl
                    B_\rmz\\
                    C\\
                    CA
                \matr
            =
                \matl
                    0 & 1 & 0\\
                    5 & 1 & 0\\
                    0 & 5 & 1
                \matr, \nn
    \end{align}
    and thus
    \begin{align}
        \SA
            &= \left[ \begin{array}{c|cc}
                    -5 & 0 & 1\\
                    \hline
                    0 & 0 & 1\\
                    0 & -6 & -5
            \end{array}\right]
        \SB=\left[
        \begin{array}{c}
            0\\
            \hline
            0\\
            1
        \end{array}
        \right], \nn
        \\
        \SC
            &=\left[    \begin{array}{c|cc}
                    0 & 1 & 0
                \end{array}
                \right].\nn
    \end{align}
    Note that $\SA_\eta = -5$ and thus ${\rm mspec}(\SA_\eta) = \{-5\},$ which confirms Theorem \ref{theorem:zeros}.
    This example shows that Theorem \ref{theorem:zeros} is unaffected by the state-space realization's lack of observability or controllability. 
    \hfill{\huge$\diamond$}
\end{exmp}

\begin{exmp}
    \label{exmp:square_MIMO_2_zeros}    
    \textbf{[Square MIMO System.]}
    Consider the square MIMO system
    \begin{align}
        G(s)=\dfrac{64}{s^3+24s^2+176s + 384}\matl
            1 & s + 4\\
            s - 2 & s-8
        \matr,\nn
    \end{align}
    with two inputs and two outputs.
    Note that $\rho_1 =2$ and $\rho_2 =2,$ and thus $\rho = 4.$
    A realization of $G$, computed with MATLAB's \texttt{ss} routine, is 
    \begin{align*}
        A
            &=
                \matl
                    -24 & -11 & -6 & 0 & 0 & 0\\
                    16 & 0 & 0 & 0 & 0 & 0\\
                    0 & 4 & 0 & 0 & 0 & 0\\
                    0 & 0 & 0 & -24 & -11 & -6\\
                    0 & 0 & 0 & 16 & 0 & 0\\
                    0 & 0 & 0 & 0 & 4 & 0
                \matr, 
        B
            =
                \matl
                    2 & 0\\
                    0 & 0\\
                    0 & 0\\
                    0 & 4\\
                    0 & 0\\
                    0 & 0
                \matr, 
        \\
        C
            &=
                \matl
                    0 & 0 & 0.5 & 0 & 1 & 1\\
                    0 & 2 & -1 & 0 & 1 & -2
                \matr.
    \end{align*}
    Furthermore, ${\rm izeros} 
        \left[
        \begin{array}{c|c}
            A & B \\
            \hline
            C & D
        \end{array} 
        \right] = \{ -1,0\},$ which is computed using the MATLAB's \texttt{tzero} routine.

    Letting $B_\rmz = \matl 0 & 0 & 1 & 0 & 0 & 0\\
    0 & -1 & 0 & 0 & 0 & 0\matr,$ it follows from \eqref{eq:T_def} that 
    \begin{align}
        T
            =
                \matl
                    B_\rmz\\
                    C_1\\
                    C_1A\\
                    C_2\\
                    C_2A
                \matr
            =
                \matl
                                0&   0&   1&   0&   0&   0\\
                                0&  -1&  0&   0&   0&   0\\
                                0&   0&   0.5& 0&   1&   1\\
                                0&   2&   0&   16&  4&   0\\
                                0&   2& 
                                -1&  0&   1&   -2\\
                                32&  -4&   0&  16&  -8&   0
                \matr,\nn
    \end{align}
    and thus
    \begin{align*}
        \SA
            &=\left[ \begin{array}{cc|cccc}
                     0&   -4&  0&   0&   0&   0\\
                     0&   -1&  -4&  0.5& -2&  -0.5\\
                     \hline
                     0&   0&   0&   1&   0&   0\\
                     48& -38& -88& -21& 4&   1\\
                     0&   0&   0&   0&   0&  1\\
                    -144&  268& -272&   -6 &  -88&   -26
            \end{array}
            \right],\nn\\
        \SB 
            &=\left[
        \begin{array}{cc}
            0 & 0\\
            0 & 0\\
            \hline
            0 & 0\\
            0 & 64\\
            0 & 0\\
            64 & 64
        \end{array}
        \right],\quad
        \SC
            =\left[    \begin{array}{cc|cccc}
                                        0&     0&     1&     0&     0&     0\\
                                        0&     0&     0&     0&     1&     0
                \end{array}
                \right].\nn
    \end{align*}
    Note that $\SA_\eta = \matl 0 & -4\\
    0 & -1
    \matr$ and thus ${\rm mspec}(\SA_\eta) = \{-1,0\},$ which confirms
    that the invariant zeros of the system are the eigenvalues of $\SA_\eta.$
    \hfill{\huge$\diamond$}
\end{exmp}

\begin{exmp}
\label{exmp:wide_MIMO}
\textbf{[Wide MIMO System.]}
    Consider the wide MIMO system
    \begin{align}
        G(s)
            =
                \dfrac{1}{d(s)} 
                \matl
                    s^2+s-2 & 0 & s^2-2s+1\\
                    s^2-3s+2 & s^2-1 & s^2-1\\
                \matr
                ,\nn
    \end{align}
    where $d(s) = s^3+2s^2-s-2,$
    with three inputs and two outputs.
    Note that $\rho_1 =1$, and $\rho_2 =1,$ and thus $\rho = 2.$
    A realization of $G$, computed with MATLAB's \texttt{ss} routine, is 
    \begin{align*}
        A
            &=
                \matl
                    0 & 0 & 2 & 0 & 0 & 0\\
                    1 & 0 & 1 & 0 & 0 & 0\\
                    0 & 1 & -2 & 0 & 0 & 0\\
                    0 & 0 & 0 & 0 & 0 & 2\\
                    0 & 0 & 0 & 1 & 0 & 1\\
                    0 & 0 & 0 & 0 & 1 & -2
                \matr, 
        \\
        B
            &=
                \matl
                    -2 & 0 & 1\\
                    1 & 0 & -2\\
                    1 & 0 & 1\\
                    2 & -1 & -1\\
                    -3 & 0 & 0\\
                    1 & 1 & 1
                \matr, 
        \\
        C
            &=
                \matl
                    0 & 0 & 1 & 0 & 0 & 0\\
                    0 & 0 & 0 & 0 & 0 & 1
                \matr.
    \end{align*}
    Note that ${\rm izeros} 
        \left[
        \begin{array}{c|c}
            A & B \\
            \hline
            C & D
        \end{array} 
        \right] = \{ 1,1\}.$     
    It follows from \eqref{eq:T_def} that 
    \begin{align}
        &T
            =
                \matl
                    B_\rmz\\
                    C_1\\
                    C_2
                \matr
            =\nn\\
                &\matl
                    0.642 &0.423 & 0.205 & 0.406 & 0.187 & 0.406\\
                     -0.532 & -0.033 & 0.465 & 0.167 & 0.665 & 0.167\\
                     -0.110 & -0.390 & -0.670 & 0.426 & 0.146 & 0.426\\
                     -0.5 & 0.5 & 0 & 0.5 & -0.5 & 0 \\
                     0 & 0 & 1 & 0 & 0 & 0\\
                    0 & 0 & 0 & 0 & 0 & 1
                \matr,\nn
    \end{align}
    where $B_\rmz$ is computed with the MATLAB routine \texttt{null(B')}.
    Finally, the zero-subspace form of $G$ is
    \begin{align*}
        &\SA
            =\\
                &\left[\begin{array}{cccc|cc}
                    0.696  &  0.265  &  -0.07 &  -0.218  &  0.984  & -0.109\\
                    0.691  &  0.394 & 0.163 &   0.498  & -2.244  &  0.249\\
                    -0.388  &  0.339  & 0.908  & -0.280 & 1.259  & -0.140\\
                    0.067 &  -0.275  &  -0.542 &  -0.5  & -0.75  & 0.75\\
                    \hline
                    0.578  & 0.506  & -1.335  &  0.5  & -3.25   & 0.25\\
                    0.442  &  1.056  & -0.249  & -0.5  & -0.75  &  -2.25
                \end{array}\right],
        \nn\\
        \SB&=\left[\begin{array}{ccc}
            0 & 0 & 0\\
            0 & 0 & 0\\
            0 & 0 & 0\\
            4 &  -0.5 &  -2\\
            \hline
            1 & 0 & 1\\
            1 & 1 & 1
        \end{array}\right],\nn 
        \\
        \SC
            &=
                \left[\begin{array}{cccc|cc}
                    0   &  0 &     0 &    0   &  1 &    0\\
                    0   &  0 &     0 &    0   &  0 &    1
                \end{array}\right].\nn
    \end{align*}
    Note that ${\rm mspec}(\SA_\eta) = \{0,-0.5,1,1\},$ which suggests
    that the invariant zeros of the system are contained in the ${\rm mspec}(\SA_\eta).$
    

    To identify the invariant zero in ${\rm mspec}(\SA_\eta),$ 
    we use the fact that the rank of the associated Rosenbrok system matrix $\SZ(s)$ drops at the invariant zero.     
    Since
    $
        \rank \SZ(0) = \rank \SZ(-0.5) = 8, 
    $
    but 
    $\rank \SZ(1) = 7,$ 
    it follows that the invariant zeros of the system are at $z=\{1,1\}$ as the rank of the Rosenbrok system matrix drops only at $z = 1.$ 
    \hfill{\huge$\diamond$}
\end{exmp}

\section{Conclusion}
\label{sec:conclusion}

This paper showed that similar to the poles computation, zeros of a state-space realization of a linear system can be computed by solving an eigenvalue problem instead of a generalized eigenvalue problem.
%
%
%
By transforming the system into its zero-subspace form, it is shown that the zeros are the eigenvalues of a partition of the dynamics matrix represented in the zero-subspace form, which is also the dynamics matrix of the system's zero dynamics.  
Numerical examples validate the main result of the paper.

Future research is focused on extending the main result of this paper to nonsquare MIMO systems. 
Although the technique is easily extended to the case of wide systems, where it provides spurious invariant zeros, which can be removed heuristically, it provides inconclusive results in tall systems. 
%
%
An alternative approach may be to consider random bordering to square the nonsquare system and thus compute invariant zeros of the MIMO system. 



\printbibliography

\section*{Appendix A: Normal Form}
\label{appndx:normalForm}
This section reviews the construction of the normal form of a nonlinear affine system \cite{kolavennu2001nonlinear,isidori1985nonlinear}.

Consider an affine system
\begin{align}
    \dot x 
        &= 
            f(x) + g(x) u,
    \label{eq:xdot_gen}
    \\
    y 
        &=
            h(x),
    \label{eq:y_gen}
\end{align}
where $x(t)\in \BBR^{l_x}$ is the state, 
$u(t)\in \BBR^{l_u}$ is the input, 
$y(t)\in \BBR^{l_y}$ is the output,
and $f, g, h$ are smooth functions of appropriate dimensions. 
For each $i\in \{ 1, \ldots, l_y \},$ let $\rho_i$ denote that relative degree of the $i$th output $y_i.$
Furthermore, let $\rho \isdef \sum_i^{l_y} \rho_i $ denote the relative degree of the system.




Consider the transformation
\begin{align}\label{eq_T}
    T : \BBR^{l_x} &\to \BBR^{l_x}
    \nn \\
        T(x)
            &=
                \matl
                    \phi(x) \\
                    \psi(x)
                \matr,
\end{align}
where $\phi(x) $ is chosen to satisfy
\begin{align}\label{eq_phi}
    L_g \phi(x) = 0,
\end{align}
and
\begin{align}
    \psi(x) 
        =
            \matl
                \psi_1(x) \\
                \vdots \\
                \psi_{l_y}(x)
            \matr,
    \label{eq:psi_def}
\end{align}
where
\begin{align}
    \psi_i(x)
        \isdef
            \matl
                h_i(x) \\
                L_f h_i(x) \\
                \vdots \\
                L_f^{\rho_i-1} h_i(x)
            \matr
        \in 
            \BBR^{\rho_i}.
    \label{eq:psi_i_def}
\end{align}
Note that $\phi : \BBR^{l_x} \to \BBR^{{l_x} - \rho}$ and, for $i=1, \ldots, l_y,$
$\psi_i : \BBR^{l_x} \to \BBR^ {\rho_i},$ and thus
$\psi : \BBR^{l_x} \to \BBR^ \rho.$
Furthermore, the functions $\psi_i$ are well-defined since the functions $f,g,h$ are assumed to be smooth. 
However, $\phi$ satisfying \eqref{eq_phi} may or may not exist. 

Assuming that $\phi$ satisfying \eqref{eq_phi} exists and defining $\eta \isdef \phi(x), $ it follows 
that 
\begin{align}
    \dot \eta 
        &=
            L_f \phi(x) + L_g \phi(x) u
        =
            L_f \phi(x),
    \label{eq_eta_dot}
\end{align}
where $L_g \phi(x)=0$ by construction. 
Note that \eqref{eq_eta_dot} is the \textit{zero dynamics} \cite{Khalil:1173048}.

Next, defining $\xi \isdef \psi(x),$ it follows
that 
\begin{align}
    \dot \xi
        &=
            L_f \psi(x) + L_g \psi(x) u.
    \label{eq:xi_dot}
\end{align}
Next, note that 
\begin{align}
    L_f \psi(x)
        =
            A_\rmc \xi 
            +
            B_\rmc 
            \matl 
                L_f^{\rho_1} h_1(x) \\
                \vdots \\
                L_f^{\rho_{l_y}} h_{l_y}(x)
            \matr,
    \label{eq:Lf_psi}
\end{align}
where $A_\rmc = {\rm diag  } (A_{\rmc,1}, \ldots, A_{\rmc, l_y}) \in \BBR^{\rho \times \rho }$ and 
$B_\rmc = {\rm diag} (b_{\rmc,1}, \ldots, b_{\rmc,l_y}) \in \BBR^{\rho \times l_y} $ and, for $i=1,\ldots, l_y,$ 
\begin{align}
    A_{\rmc,i}
        &\isdef 
            \matl 
                0 & 1 & 0 & \cdots & 0  \\
                0 & 0 & 1 & \cdots & 0  \\
                \vdots & \vdots & \ddots & \ddots & \vdots \\
                0 & \vdots & \ldots & 0 & 1  \\
                0 & \vdots & \ldots & 0 & 0  \\
            \matr
        \in \BBR^{\rho_i \times \rho_i},
    \\
    b_{\rmc,i}
        &\isdef 
            \matl 
                0\\
                \vdots \\
                1
            \matr
        \in \BBR^{\rho_i}.
\end{align}
Finally, 
\begin{align}
    L_g \psi(x)
        &=
            B_\rmc  
            \matl 
                L_g L_f^{\rho_1-1} h_1(x) \\
                \vdots \\
                L_g L_f^{\rho_{l_y}-1} h_{l_y} (x)
            \matr.
    \label{eq:Lg_psi}
\end{align}

Substituting \eqref{eq:Lf_psi} and \eqref{eq:Lg_psi} in \eqref{eq:xi_dot} yields
\begin{align}
    \dot \xi
        &=
            A_\rmc \xi 
            +
            B_\rmc
            \left(
                \alpha(x)
                +
                \beta(x) u
            \right),
\end{align}
where 
\begin{align}
    \alpha(x)
        &\isdef
            \matl 
                L_f^{\rho_1} h_1(x) \\
                \vdots \\
                L_f^{\rho_{l_y}} h_{l_y}(x)
            \matr
            \in \BBR^{l_y}, 
    \label{eq:alpha_def}
    \\
    \beta(x)
        &\isdef 
            \matl 
                L_g L_f^{\rho_1-1} h_1(x) \\
                \vdots \\
                L_g L_f^{\rho_{l_y}-1} h_{l_y} (x)
            \matr
            \in \BBR^{l_y \times l_u}.
    \label{eq:beta_def}
\end{align}

The normal form of the nonlinear system \eqref{eq:xdot_gen}, \eqref{eq:y_gen} is then 
\begin{align}
    \dot \eta 
        &=
            L_f \phi(x),
        \\
    \dot \xi
        &=
            A_\rmc \xi 
            +
            B_\rmc
            \left(
                \alpha(x)
                +
                \beta(x) u
            \right).
\end{align}


\end{document}